\documentclass[11pt]{article}
\usepackage{epsfig}
\usepackage{color}
\usepackage{amsmath}
\usepackage{amsthm}
\usepackage{amssymb}

\usepackage{amscd}
\usepackage{amsfonts}  
\usepackage{hyperref}
\usepackage[all]{xy}
\usepackage{graphicx}    
\usepackage{multicol}    

\DeclareMathAlphabet\EuR{U}{eur}{m}{n}
\SetMathAlphabet\EuR{bold}{U}{eur}{b}{n}

\DeclareMathOperator{\im}{im}

\newcommand{\boundary}{\partial}
\newcounter{commentcounter}
\newcommand{\comment}[1]                      
{\stepcounter{commentcounter}
{\bf Comment \arabic{commentcounter}}: {\ttfamily #1} }

\newcommand{\bbR}{{\mathbb R}}

\newcommand{\bbZ}{{\mathbb Z}}

\newcommand{\id}{\operatorname{id}}

\newcommand{\st}{\mathbf{S}}

\theoremstyle{plain}
\newtheorem{theorem}{Theorem}
\newtheorem{lemma}[theorem]{Lemma}
\newtheorem{proposition}[theorem]{Proposition}

\theoremstyle{definition}
\newtheorem{definition}[theorem]{Definition}

\theoremstyle{remark}
\newtheorem{remark}[theorem]{Remark}

{\catcode`@=11\global\let\c@equation=\c@theorem}

\title {A geometric description of smooth cohomology
}
\author{Ulrich Bunke, Matthias Kreck and Thomas Schick}

\begin{document}
\maketitle

In this paper we give a geometric cobordism description of smooth integral
cohomology. The main motivation to consider this model (for other models see
\cite{cheeger85:_differ}, \cite{MR2192936}, \cite{Dupont}) is that it allows
for simple descriptions of both the cup product and the integration, so that
it is easy to verify
the compatibilty of these structures. We proceed  in a similar way in the case
of  smooth cobordism as constructed in \cite{B-S-S-W}. There the starting
point was Quillen's cobordism description of singular cobordism groups for a
smooth manifold $X$. Here we use instead the similar description of ordinary
cohomology from \cite {K}. This cohomology theory is denoted by $SH^k(X)$. In
this description smooth manifolds in Quillens' description are replaced by
so-called stratifolds, which are certain stratified spaces. The cohomology theory $SH^k(X)$ is naturally
isomorphic to ordinary cohomology $H^k(X)$, thus we obtain a cobordism type
definition of  the smooth extension of ordinary integral cohomology.

\section{Axioms of smooth cohomology theories}

To begin, let us recall what is meant by a \emph{smooth extension of
 the functor  $H^k$}, ordinary cohomology. Compare \cite[Definition 1.1]{bunke07:_smoot_k_theor} for the
formal definition and the
fundamental paper \cite{MR2192936} for a general construction of such
theories. We denote the closed $k$-forms by $\Omega_{cl}^k(X)$ and the map to
de Rham cohomology, which we identify with real singular cohomology via the de
Rham isomorphism, by
$Rham\colon  \Omega_{cl}^k(X)\to H^k(X;\mathbb R) $.\\

\begin{definition}\label{def:smoothext}
 A \emph{ smooth extension} of $H$ is a
    functor $X \mapsto \hat{H}^*(X)$ from the category of smooth manifolds to
    $\mathbb Z$-graded groups together with natural transformations
  \begin{enumerate} \item $R \colon
      \hat{H}^*(X) \to \Omega_{cl} ^*(X)$ (curvature) \item $I \colon
      \hat{H}^*(X) \to H^*(X)$ (forget smooth data) \item $a\colon \Omega
      ^{*-1}(X)/{\im(d)}\to \hat {H}^*(X)$ (action of forms). 
      \end{enumerate}

These transformations have to
  satisfy the following axioms.
  \begin{enumerate}
  \item  $R\circ a = d\colon \Omega^{*-1}(X) \to \Omega_{cl}^*(X)$.
\item  The following diagram commutes:
$$\xymatrix{\hat {H}^*(X)\ar[d]^R\ar[r]^I&H^*(X)\ar[d]\\\Omega_{cl}^k(X)\ar[r]^{Rham}&H^k(X;\mathbb R)}.$$

  \item  For every smooth manifold $X$ the sequence
$$H^{*-1}(X) \to \Omega ^{*-1}(X)/{\im(d)}\xrightarrow{a} \hat {H}^*(X) \xrightarrow{I} H^*(X) \to 0$$
is exact, where the first map is the composition 
$$H^{*-1}(X)\longrightarrow  \frac{ker(d\colon \Omega^{*-1}(X)
\to\Omega^*(X))}{\im (d)} \subseteq \frac{\Omega^{*-1}(X)}{\im(d)}\ .$$
\end{enumerate}

To \emph{have a compatible ring structure} means that $\hat{H}$ actually
takes values in graded commutative rings (we denote the product by $\cup$), that $R,I$ are ring maps, and that
for all $x \in \hat{H}(X)$  and $\omega \in \Omega(X)/{\im(  d)}$ we have 
$$
a(\omega) \cup x = a(\omega \wedge R(x)).
$$
In this case we call the smooth extension a \emph{multiplicative smooth
  extension}. 
  
  We will use a construction of $H^k(X)$ in terms of cobordism classes similar to Quillen's description of singular cobordism. The difference is that we replace manifolds by manifolds with singularities called stratifolds. In the next section we briefly introduce stratifolds and prove some basic properties needed for our construction. To distinguish the cohomology theory constructed from stratifolds from singular cohomology we denote it by $SH^k(X)$.

\end{definition}

\section{Stratifolds}

Now we give a short introduction to stratifolds, where we are rather sketchy
and refer the reader to \cite{K} for details.  A stratifold is a topological
space $\st$ together with a subsheaf $\mathcal C$ of the sheaf of continuous
functions, which in case of a smooth manifold plays the role of the sheaf of
smooth functions. The space $\st$ and the sheaf $\mathcal C$ have to fulfill
certain natural axioms, which in particular give a decomposition of $\st$ into
smooth manifolds, the \emph{strata} of $\st$. The top-dimensional stratum
$\st^m$ is also called the regular part $\st^{reg}$.
There is an obvious definition of \emph{morphisms} between stratifolds, which
are continuous maps $f\colon \st \to \st'$, which pull elements in $\mathcal
C'$ back to elements in $\mathcal C$.  A basic property which relates the
strata $\st^k$ to the stratifold $\st$ is that for each $x\in \st^k$ there is
an open neighborhood $U\subseteq \st$ of $x$  and a retract $r\colon U \to U \cap
\st^k$ which is a morphism.  These are called local retracts. It is useful to
note that a map is a morphism if and only if it is smooth on all strata and commutes with appropriate local retracts.
We will only consider \emph{regular stratifolds} which means that locally near
each $x $ in a stratum $\st^k$  the stratifold looks like $V \times F$ for
some stratifold $F$ and open subset $V$ in $\st^k$. A stratifold is called
\emph{orientable} if the codimension $1$ stratum is empty and the top stratum
is orientable. Once an orientation on the top stratum is fixed we call such a
stratifold an \emph{oriented stratifold}. More generally, if $\st$ is a
stratifold and $X$ a smooth manifold, a continuous map $f\colon \st\to X$ is
called orientable if the codimension $1$ stratum of $\st$ is empty and if
$f|_{\st^{reg}}^*\Lambda_X^{\bbZ/2}$ is isomorphic to
$\Lambda_{\st^{reg}}^{\bbZ/2}$. Here, $\Lambda_X^{\bbZ/2}$ is the orientation
covering of $X$. An \emph{orientation of $f$} is then the choice of such an
isomorphism. Note that it also gives an isomorphism between
$\Lambda_{\st^{reg}}$ and $f|_{\st^{ref}}^*\Lambda_X$, where $\Lambda_X$ is
the real orientation bundle of $X$.

 We also consider stratifolds with boundary. This is a pair of spaces  $(\mathcal T, \partial \mathcal T)$ together with stratifold structures on $\mathcal T - \partial \mathcal T$ and on $\partial \mathcal T$ and a germ of collars $c \colon  \partial \mathcal T \times [0,1) \to \mathcal T$. 
Many basic properties of smooth manifolds generalize to stratifolds, like tangent spaces (the vector space of derivations of the germ of morphisms to $\mathbb R$), the differential of a morphism, differential forms (see below), Sard's theorem, approximation of continuous maps from a stratifold to a smooth manifold by morphisms and the transversality theorem for a map from a stratifold $\st$ to a smooth manifold $X$ and a smooth map from a manifold $Y \to X$. For all this see \cite {K}.\\

Since we will use differential forms intensively we define them on stratifolds. A $k$-form $\omega$  on a stratifold $\st$ is a prescription which to each $x \in \st$ assigns an alternating $k$-form on $T_x\st$, which fulfills the following property:
\begin{enumerate}
\item 
 The restriction to each stratum is a differential $k$-form.
\item For each $x \in \st^r$, the $r$-stratum, there is an open neighborhood $U\subseteq \st$
  of $x$ in $\st$ and a local retract
$$
r\colon U \to U \cap \st^{m},
$$
such that
$$
\omega|_{U} = r^*(\omega|_{U \cap \st^r}).
$$
\end{enumerate}
Here the pull back of a differential form under a morphism from a stratifold
to a smooth manifold is defined as for smooth manifolds using the
differential.

\begin{lemma}\label{lem:support}
  Let $\st$ be an $m$-dimensional stratifold and $f \colon  \st
  \to X $ be a proper morphism  to an $n$-dimensional  smooth manifold.  Then there is an open
    neighborhood $V$ of the singular part such that for all $\omega \in
    \Omega^{m} _c(X;E)$ the pullback $f^*(\omega) $
vanishes on $V$. \textcolor{black}{Here, $E\to X$ is any coefficient bundle. In other words, there is a fixed compact subset
$K=\st\setminus V$ of the regular part such that $supp(f^*\omega)\subset K$.}
\end{lemma}

\begin{proof} Let $x\in \st^r$ for $r< m$ be a point in the singular part of $\st$. Then there is an open neighborhood of $x$ in $\st$ and a local retract 
$$
r\colon 
U \to U \cap \st^{r},
$$
such that the restriction of $f$ to $U$ factors over $r$ and hence
$$
f^*(\omega|_{U}) = r^*(f^*(\omega)|_{U \cap \st^r}).
$$
Note that $U$ and $r$ are determined by $f$ and can be chosen independent of
$\omega$. 

Since $f^*(\omega)|_{U \cap \st^r} =0$ for dimensional reasons we conclude that there is an open neighborhood $V$ of the singular part $\st - \st^{m}$ such that $
(f|_{\st ^{m}})^*(\omega) 
$ vanishes on $\st^{m} \cap V$. \\

Now, let $K$ be a compact set in $X$ such that $\omega$ vanishes outside of $K$.  Then, since $f$ is proper, $f^{-1}(K)$ is a compact subset of $\st$. Since $\st ^{m} - V$ is a closed subset of $\st^{m}$, the set $f^{-1}(K)\cap (\st^{m} - V)$ is compact, and $f^*(\omega)$ vanishes outside this set.
\end{proof}

As a consequence we can define the following integral. Let $\omega$ be a form
with compact support on  $X$ \textcolor{black}{with coefficients in $\Lambda_X$} and $f \colon  \st \to X$ a proper \textcolor{black}{oriented}
morphism from an $m$-dimensional stratifold $\st$ to $X$. Then we define
$$
\int _\st f^*(\omega) := \int_{\st^{m}}
(f|_{\st^m})^*(\omega),
$$
\textcolor{black}{using the identification of $f|_{\st^m}^*\Lambda_X$ with
  $\Lambda_{\st^m}$ from the orientation of $f$.}

Similarly, for a proper \textcolor{black}{oriented} morphism $F\colon  \mathcal{T} \to X$ from an
$m$-dimensional stratifold with boundary we define
$$
\int _\mathcal T F^*(\omega) := \int_{\mathcal T^{m}}
(F|_{\textcolor{black}{\mathcal T}^m})^*(\omega).
$$ 
Stokes' Theorem applied to the top stratum gives us {\bf Stokes' Theorem}: If $\omega \in \Omega^{m-1}(X;\Lambda_X)$, then
\begin{equation}\label{2theorem}
\int_{\mathcal T}d(F^*(\omega)) = \int _{\partial \mathcal T}F|_{\boundary\mathcal{T}}^*(\omega).
\end{equation}

\begin{remark}
  Throughout this article, we will work with oriented morphisms and maps
  $f\colon \st\to X$ and with differential forms with values in the
  orientation bundle of $X$. In the special case that $X$ is an oriented
  manifold, these are ordinary forms, and an orientation of $f$ is precisely
  an orientation of $\st$. A reader not used to the more general setting might
  just assume the orientability and choice of orientations throughout. The
  passage to the general case is a mere technical point.
\end{remark}
\section{Differential cohomology via stratifolds}

Now we define the cycles of our smooth cohomology following the \textcolor{black}{recipe}
as for for singular cobordism \cite {B-S-S-W}. The starting point
is the description of $k$-th ordinary cohomology of $X$ as bordism classes of
 continuous  oriented proper maps 
from oriented regular stratifolds $\st $ of dimension $n-k$ to $X$
\cite[Chapter 12]{K}. \textcolor{black}{Actually, in \cite{K} one has to make the obvious
modifications to pass from oriented manifolds to arbitrary manifolds by
working with oriented maps.}  We call such \textcolor{black}{oriented} proper maps \emph{$k$-cycles}. We denote this
cohomology group $SH^k(X)$, the \emph{stratifold cohomology} of $X$. 
Since every proper map is homotopic
via proper maps to a morphism \cite{K} we will always assume that $f$ is a
morphism.  

Let $f\colon  \st \to X$ be a $k$-cycle, i.e.~a proper
\textcolor{black}{oriented} morphism from a
regular $(n-k)$-dimensional stratifold to $X$. Then we construct a current
$
T(\st,f)
$, i.e. an element
in the topological dual space ${\Omega}^{n-k}_c(X;\Lambda_X) ^*$ of continuous
linear maps from $ {\Omega}^{n-k}_c(X;\Lambda_X)$ to $\mathbb R$ as follows:
$$
\omega \mapsto \int_{\st} f^*(\omega).
$$
\\
We have an injective map $$
j\colon   \Omega ^k(X) \to  \Omega^{n-k}_c(X;\Lambda_X) ^*
$$
given by
$$j(\alpha):=\{\omega\mapsto \int_X \alpha\wedge \omega\}\ .$$

After these preparations we define a cycle for $k$-th smooth cohomology:\\

\begin{definition}
  A {\bf smooth cycle of degree $k$} is a triple
$$
(\st,f,\alpha)
$$
where $(\st,f)$ is as above given by a proper \textcolor{black}{oriented} morphism $f\colon \st \to X$,
with $\st$  an $(n-k)$-dimensional oriented regular stratifold, and  $\alpha
\in  \Omega^{n-(k+1)}_c(X;\textcolor{black}{\Lambda_X})^*/{\im(d^*)}$, such that $T(\st,f) - d^*(\alpha)$
is in the image of $j$. \textcolor{black}{The sum of two smooth cycles is
  defined by disjoint union. The negative of a cycle $(\st,f,\alpha)$ is
  $(\st,f^-,-\alpha)$ where $f^-$ is $f$ with the reverse orientation.}

We define cobordisms of cycles for smooth cohomology as follows: if $\mathcal{T}$ is a stratifold with boundary $\boundary 
\mathcal{T}=\st$ and $F\colon \mathcal{T}\to X$ is a proper
\textcolor{black}{oriented} morphism, we say that $(\mathcal{T},F)$ is a zero
bordism of $(\boundary W,F|_{\boundary \mathcal{T}},T(\mathcal{T,F}))$, and a
bordism between two cycles is a zero bordism of the difference. The 
only thing which is a bit special is that for the map $F$  on a bordism
$\mathcal T$ to $X$ we require that $F$ commutes with the collar $c\colon
\partial T \times [0,\epsilon)$ for some appropriate $\epsilon >0$. This
allows to glue bordisms (using the collars) in a compatible way. Note also
that $dT(\mathcal{T,F})-T(\boundary \mathcal{T},F|_{\boundary \mathcal T})=0$
(in particular is {in the image of $j$}) by
(\ref{2theorem}): 
$$d^{{*}}T(\mathcal{T,F})(\omega)=\int_{\mathcal{T}}F^*(d\omega) = \int_{\mathcal
  T}dF^*\omega 
=\int_{\boundary\mathcal{T}}F|_{\boundary \mathcal
  T}^*\omega=T(\boundary\mathcal{T},F|_{\boundary\mathcal T})(\omega)\ .$$ 

We denote the corresponding bordism group by
$$
\hat {SH}^k (X).
$$
We call this group the \emph{smooth stratifold cohomology of $X$}.
 
As in \cite[\textcolor{black}{Definition 4.9}]{B-S-S-W} we define the maps 
  \begin{eqnarray*}
    R \colon \hat{SH}^k(X) \to \Omega^k(X); &
    (\st,f,\alpha)\mapsto j^{-1}(T(\st,f)-d^*(\alpha)),\\
    a\colon \Omega^{k-1}(X)/{\im(d)}\to \hat {SH}^k(X); &
    \alpha \mapsto (\emptyset,-j(\alpha)),\\
    I\colon \hat {SH}^k(X) \to SH^k(X);&
 (\st,f,\alpha)\mapsto[\st,\alpha].
\end{eqnarray*}
\end{definition}

The proof that these maps are well defined is literally the same as in the
case, where we have smooth manifolds instead of stratifolds \cite[\textcolor{black}{Lemma 4.10}]{B-S-S-W},
since the basic ingredient, Stokes' Theorem, is available.



The next aim is to construct induced maps for a smooth map $g\colon  Y \to
X$. The basic idea is, that if $(\st,f,\alpha)$ is a smooth cycle in $X$, then
we can after a proper homotopy of $f$ (which we can consider as a special case of a bordism)
assume that $g$ is transversal to $f$.
Then, as in \cite{K}, one can consider the pull back of $\st$ giving a cycle
$(g^*(\st),F)$ in $Y$. We denote the canonical map $g^*(\st) \to \st$ by
$G$. \textcolor{black}{As in \cite[\textcolor{black}{Section 4.2.6}]{B-S-S-W},
  the orientation of $f$ induces an orientation of $F=g^*f$.}

To extend this pull back to a smooth cycle by pulling back $\alpha$, one
has the same situation as in \cite{B-S-S-W}, i.e.~one has to pull back
$\alpha$ along $g$. Recall that this is only possible if $WF(\alpha)\cap
N(g)=\emptyset$. Here $WF(\alpha)\subseteq T^*X$ denotes the wave front set of
the distributional form $\alpha$, and $N(f)\subseteq T^*X$ is the normal set
to $f$. \textcolor{black}{The
  wave front set of  a distributional form $\alpha$  on $X$ is a conical
  subset of $T^*X$ which measures the locus and the directions of the
  singularities of $\alpha$. For a precise definition and for the properties
  of distributions using the wave front set needed we refer to
  \cite[{Section 8}]{MR1996773}.} Compare \cite[Section 4.2.6]{B-S-S-W} for the 
notation and more 
details. In terms of normal sets transversality of  $f$ and $g$ can be expressed as
$N(f)\cap N(g)=\emptyset$ (where $N(f)$ is the normal set of the restriction of $f$ to the top stratum). 
Hence $g^*\alpha$ is defined if 
$WF(\alpha)\subseteq N(f)$. In order to match this condition we use the freedom to change $\alpha$ by an element in
the image of $d^*$.

We observe that by definition and \textcolor{black}{by Lemma \ref{lem:support}} $T(\st,f)=(f|_\st^{reg})_!(\rho)$ where
$\rho \colon \st^{reg}\to \bbR$ is a smooth \textcolor{black}{compactly supported} cutoff function which is zero in
a neighborhood of the singular set of $\st$, and which is identically $1$ on
the support of any $f^*\omega$ for $\omega\in \Omega_c^{n-k}(X)$.

Consequently, the construction of $T$ is described entirely in the context of
smooth manifolds, smooth maps and smooth forms; as in the context of
\cite{B-S-S-W}. Now, the arguments there, in particular \cite[Lemma
4.12]{B-S-S-W} literally apply in
our situation to show that we can change $\alpha$ to a representative $\alpha
'$ 
such that the wave front  \textcolor{black}{set} of $\alpha$ satisfies $WF(\alpha') \subseteq
N(f)$. Then $g^*(\alpha')$ is a well defined distribution and we can make the
following definition.

\begin{definition}
  We set $g^*[\st,f,\alpha] :=[g^*(\st),\textcolor{black}{g^*f}, g^*(\alpha)]$, where we choose a
  representative such  that
   \textcolor{black}{$f$ is transversal to $g$} and $WF(\alpha) \subseteq N(f)$.
\end{definition}

The proof that this induced map is well defined and functorial is the same as
in the case where $\st$ is a smooth manifold. \textcolor{black}{Naturality of the transformations $R$, $I$ and
$a$ is checked in a straigtforward way}.

\section{Ring structure on smooth \protect\textcolor{black}{stratifold cohomology}}

\begin{definition}
  We define the \emph{$\times$-product} of classes $[\st,f,\alpha] \in
  \hat{SH}^k(X)$ and $[\st',f',\alpha'] \in \hat{SH}^r(X')$ \textcolor{black}{with values in
  $\hat{SH}^{k+r}(X\times X')$} as
  \begin{multline*}
    [\st,f,\alpha]\times [\st',f',\alpha'] :=
    [(-1)^{kr}\st \times \st', f\times
f', (-1)^kR([\st,f,\alpha]) \times \alpha' + \alpha \times T(\st',f')]
.
\end{multline*}
The sign $(-1)^{kr}$ comes from the fact that in contrast to \cite{B-S-S-W} we
work with orientations of the tangent bundle, whereas there
normal orientations \textcolor{black}{are used}. This orientation convention
is in agreement with that in \cite{K}.
\end{definition}
The proof of the following fundamental properties is the same as in \cite{B-S-S-W}, except that for the difference of signs one has to use the arguments in \cite{K}.

\begin{proposition} The product is well defined, associative, graded commutative, and natural.
\end{proposition}

Using the map induced by the diagonal $\Delta\colon X\to X\times X$ we define
the cup product.

\begin{definition}
  For $a\in \hat{SH}^k(X)$ and $b \in \hat{SH}^r(X)$ we define
$$
a \cup b:= \Delta ^*(a \times b).
$$
\end{definition}

By a straightforward calculation we see

\begin{proposition} The maps $R$ and $I$ are multiplicative and 
$$
a(\beta) \cup [\st,f,\alpha] = a(\beta  \wedge R([\st,f,\alpha] )).
$$
\end{proposition}

\section{\textcolor{black}{Smooth stratifold cohomology as smooth extension of {ordinary} cohomology}}



\textcolor{black}{We have constructed a smooth extension of the $SH$-homology theory as
developped in \cite{K}. However, we argue that $\hat{SH}$ is a smooth exension
of ordinary cohomology.}
For this, we have to observe that the corresponding stratifold cohomology
$SH(X)$ {is naturally isomorphic to ordinary cohomology $H(X)$}. The reason is
that this functor fulfills the homotopy axiom (obvious) and that one has a
natural Mayer-Vietoris sequence. This was proven in \cite{K} for the case
where $X$ is oriented, but the same proof works in the non-oriented case. Now
we apply \cite{K-S} or \cite[Section 7]{BunkeSchickUniqueness}. There, it is proven that a cohomology theory on smooth
manifolds \textcolor{black}{is naturally isomorphic to ordinary cohomology if
it satisfies the homotopy axiom, the Mayer-Vietoris sequence and if the
cohomology groups of a $0$-dimensional manifold $X$ are the direct product 
of the cohomology groups of the points in $X$}. Moreover, it was shown in \cite{K-S} that the natural
isomorphism can be chosen to preserve the ring structure. Thus we can identify
the multiplicative cohomology theories stratifold cohomology and ordinary
cohomology.


We now formulate our main theorem

\begin{theorem}
  Our construction $\hat{SH}$ defines a multiplicative smooth extension of
  ordinary cohomology with integer coefficients in the sense of Definition
  \ref{def:smoothext}. By  \cite[Theorem 1.1]{MR2365651} or
  \cite{BunkeSchickUniqueness} it follows that
  our theory is uniquely naturally equivalent to any other of the many models for this
  extension, in particular to Cheeger-Simons differential characters of
  \cite{cheeger85:_differ}, as described in \cite{MR2365651}.

  This is actually even true as a multiplicative extension: by \cite[Theorem
  1.2]{MR2365651} or \cite{BunkeSchickUniqueness},  there is only one multiplicative smooth extension of
  ordinary cohomology. 
\end{theorem}
\begin{proof}
\textcolor{black}{Our setup is not quite
  identical to the one of \cite{MR2365651}, as there it is required that the
  kernel of $R\colon \hat{HS}^*(M)\to \Omega^*(M)$ is naturally identified
  with $H^{*-1}(M;\bbR/\bbZ)$.}

\textcolor{black}{ As this is not the case, we use instead the method of proof of
  \cite{BunkeSchickUniqueness}. There, a natural transformation $\Phi$ between any
  two smooth extensions $\hat{H}$ and $\hat{H}'$ of integral cohomology is
  constructed by making a universal
  choice. It is shown that $\Phi$ is additive and unique in even degrees because
  $H^k(pt)=0$ for $k$ odd. The same method implies immeditately that the
  transformation is additive and unique in all degree except for $*=1$, as
  $H^{*-1}(pt)=0$ except if $*=1$.}

\textcolor{black}{Next, the method shows that for the transformation $\Phi$ there is
  \begin{equation*}
  c\in\bbR/\bbZ=H^0(S^1\times S^1;\bbR/\bbZ)=H^0(K(\bbZ,1)\times K(\bbZ,1);\bbR/\bbZ)
\end{equation*}
such that for two
  classes $x,y$ of degree $1$ we have
  \begin{equation*}
  \Phi(x+y)=\Phi(x)+\Phi(y)+a(c).
\end{equation*}
However, we can modify $\Phi$ by setting
  $\Phi'(x):=\Phi(x)-a(c)$ if $x$ is of degree $1$ and $\Phi'(x)=\Phi(x)$
  otherwise. Then we conclude that $\Phi'$ is the unique additive
  transformation between the two smooth extensions of integral cohomology
  satisfying our axioms.}

\textcolor{black}{The methods of \cite{BunkeSchickUniqueness} finally show that there is at
  most one ring structure on a smooth extension of integral cohomology. Again,
  this follows because of the vanishing of $H^*(pt;\bbR/\mathbb{Z})$ if $*\ne 0$,
  together with the consideration of distributivity for products of classes of
  degree zero and degree one.}
\end{proof}

{\color{black}
\begin{proposition}
  The flat theory corresponding to $\hat{SH}$,
i.e.~the functor $U^*(X):=\ker(R\colon \hat{SH}^*(X)\to \Omega^*(X))$ is
naturally isomorphic to $H^{*-1}(X;\bbR/\bbZ)$.

In particular, $\hat{SH}$ satisfies the setup of \cite{MR2365651}.
\end{proposition}
\begin{proof}
  This is a special case of \cite[Theorem 7.12]{BunkeSchickUniqueness}. 
\end{proof}
}
\section{Integration along the fiber}

Let $p\colon X\to B$ be the projection map of a locally trivial fiber
bundle. To define ``integration along the fibers'' of $p$ for a cohomology
theory ${E}$, one has to choose an $E$-orientation for $p$ (which might not
exist).

For a general cohomology theory $E$ and a smooth extension $\hat E$, usually
one has to choose further data in addition to an ordinary $E$-orientation to
prescribe an $\hat E$-orientation, compare
e.g.~\textcolor{black}{\cite[Section 3.1]{bunke07:_smoot_k_theor} or \cite[Section 4.3.7]{B-S-S-W}}.

An exception is ordinary (integral) cohomology $H$, as already observed in
\cite{cheeger85:_differ, Dupont, Koehler}. Here, an ordinary orientation
determines canonically a smooth orientation and a smooth integration map. In
our model of smooth cohomology using stratifolds, the definition of the
integration map as well as the proof of its main properties is particularly
simple. More precisely, we will prove the following theorem. 

\begin{theorem}\label{integration}
  Given a locally trivial smooth fiber bundle with closed $d$-dimensional fibers which is oriented for ordinary cohomology, there is a canonical integration for smooth stratifold cohomology 
  \begin{equation*}
  \hat{p}_!\colon \hat{SH}^k(E)\to \hat{SH}^{k-d}(B).
  \end{equation*}
  This has the following properties
  \begin{enumerate}
  \item\label{compatibility}
  The smooth integration is compatible with integration of forms and of ordinary cohomology classes, i.e.~the following diagrams commute:
  \begin{equation*}
  \begin{CD}
  \Omega^{k-1}(E) @>{\alpha}>> \hat{SH}^k(E) @>{I}>> H^k(E)\\
    @VV{\int_{E/B}}V  @VV{\hat{p}_!}V @VV{p_!}V\\
   \Omega^{k-1-d}(B) @>{\alpha}>> \hat{SH}^{k-d}(B) @>{I}>> H^{k-d}(B).
  \end{CD}
  \end{equation*}
  \begin{equation*}
  \begin{CD}
  \hat{SH}^k(E) @>{R}>> \Omega^k_{cl}(E)\\
   @VV{\hat{p}_!}V @VV{\int_{E/B}}V\\
   \hat{SH}^{k-d}(B) @>R>> \Omega^{k-d}_{cl}(B).
  \end{CD}
  \end{equation*}
  \item\label{naturality} Naturality: If
  \begin{equation*}
  \begin{CD}
  F @>v>> E\\
  @VVqV @VVpV\\
  C@>u>> B
  \end{CD}
  \end{equation*}
  is a cartesian diagram,
   then
   \begin{equation*}
   \begin{CD}
   \hat{SH}^k(E) @>{\hat{v}^*}>> \hat{SH}^k(F)\\
   @VV{\hat{p}_!}V @VV{\hat{q}_!}V\\
   \hat{SH}^{k-d}(B) @>{\hat{u}^*}>> \hat{SH}^{k-d}(C)
   \end{CD}
   \end{equation*}
   commutes, where we use on $q$ the pullback of the orientation on $p$.
   \item\label{functorality} Functoriality: if $r\colon D\xrightarrow{q} E\xrightarrow{p} B$ is a composition of two smooth oriented fiber bundles (with composed orientation), then
   \begin{equation*}
   \hat{r}_!=\hat{p}_!\circ \hat{q}_!.
   \end{equation*}
   \item\label{projformula} Projection formula: if $x\in \hat{SH}^k(B)$, $y\in\hat{SH}^m(E)$ then
   \begin{equation*}
   \hat{p}_!(\hat{p}^*(x)\cup y) = (-1)^{kd}x\cup \hat{p}_!(y).
   \end{equation*}
   \item\label{realization} On a cycle $x=[\st,f,\alpha]$,
   \begin{equation*}
   \hat{p}_!(x) = [\st,p\circ f, \int_{E/B}\alpha:=p_*\alpha],
   \end{equation*}
   where by definition $\int_{E/B}\alpha(\omega)= \alpha(p^*\omega)$ \textcolor{black}{and where
   we equip $p\circ f$ with the composed orientation}. 
 
  \end{enumerate}
\end{theorem}


In the remainder of this section we prove Theorem \ref{integration}. 
   Note that \ref{realization} actually is a definition of $\hat{p}_!$ which
   by construction is compatible with the addition in $\hat{SH}$. However, we
   have to check that it is well defined. By compatibility with addition, for
   a cycle $(\st,f,\alpha)$ representing zero, i.e.~$\st=\boundary
   \mathcal{T}$, $f=(F\colon \mathcal{T}\to X)|_{\st}$,
   $\alpha=T(\mathcal{T},F)$, we have to check that this is mapped to zero
   under $\hat{p}_!$. However, $[\st,p\circ f,p_*\alpha]$ is precisely the
   boundary of $[\mathcal{T},p\circ F,T(\mathcal{T},F)]$, as by definition $T$ is the
   pushdown of the fundamental class, and this is natural for composition, so
   \ref{realization} indeed defines $\hat{p}_!$.
   

   Next we prove the compatibilities of \ref{compatibility}. Here, we have for a class of cycles $x=[\st,f,\alpha]$ and a form $\omega\in\Omega^*(E)$:
   \begin{equation*}
      p_! I(x)=p_![f\colon \st\to E]=[p\circ f\colon \st\to {B}]=I( \hat{p}_!(x)),
   \end{equation*}
   \begin{equation*}
   \int_{E/B} R(x) = \int_{E/B}(T(\st,f)-d^*\alpha)= T(\st,p\circ f)-d^*\int_{E/B}\alpha = R(\hat{p}_!(x)) ,
   \end{equation*}
   \begin{equation*}
   \hat{p}_!\alpha(\omega)=\hat{p}_![\emptyset,-\omega]=[\emptyset,-\int_{E/B}\omega]=\alpha(\int_{E/B}\omega). 
   \end{equation*}
   
   To prove naturality \ref{naturality} with respect to pullback in a diagram
   \begin{equation*}
   \begin{CD}
   F@>v>> E\\
   @VVqV @VVpV\\
   C@>u>> B
   \end{CD}
   \end{equation*}
   choose (without loss of generality) the cycle $x=[\st,f,\alpha]$ such that $f$ is transversal to $v$ and $WF(\alpha)\subseteq N(f)$. 
    Since $p$ and $q$ are submersions, the composition $p\circ f$ is transversal to $u$. Moreover (with the notation $WF_y(\beta)=WF(\beta)\cap T_y^*B$ and similar for normal sets) 
   $$WF_y(\int_{E/B}\alpha)\subseteq \bigcup_{x\in p^{-1}(y)} (dp_x^*)^{-1}(WF_x(\alpha))\subseteq   \bigcup_{x\in p^{-1}(y)} (dp_x^*)^{-1}N_x(f)\subseteq  N_y(p\circ f),$$
    so that $u^*(   \hat{p}_!(x))$
    is defined using the cycle $(\st,p\circ f,\int_{E/B}\alpha)$. 
   
   Then,
   \begin{multline*}
   \hat{q}_!(v^*(x)) = \hat{q}_!(v^*\st,v^*f,v^*\alpha)=(v^*\st,q\circ v^*f,\int_{F/C}v^*\alpha)\\
    = (u^*\st,u^*(p\circ f), u^*\int_{E/B}\alpha)=u^*(\hat{p}_!(x)). 
   \end{multline*}
   Here we use that  pullback and pushdown of distributional forms in a cartesion square are compatible (which follows from the corresponding statement for smooth forms by continuity, as the pullback is extended from smooth to distributional forms by continuity).
   
   It remains to prove the projection formula \ref{projformula}. This we do in two steps. First we consider the projection $\id\times p\colon B\times E\to B\times B$. If $x=[\st,f,\alpha]\in \hat{SH}^k(B)$ and $y=[\tilde\st,\tilde f,\tilde \alpha]\in \hat{SH}^m(E)$ then
   \begin{multline*}
\hspace{3cm}     \widehat{\id\times p}_!(x\times y)=\\
\
     (-1)^{km}[\st\times\tilde\st,f\times(p\circ\tilde f),(-1)^k R(S,f,\alpha)\times\int_{E/B}\tilde \alpha+\alpha\times T(\tilde S,p\circ\tilde f)]\\
    \hspace{-2cm}= (-1)^{kd} x\times \hat{p}_!(y)\ . 
   \end{multline*}
Secondly, using the diagonal inclusion $B\to B\times B$ we pull back the whole situation to $p\colon E\to B$ and use the naturality of the smooth integration with respect to pullback. Observe  that the natural map $E\xrightarrow{p\times \id} B\times E$ which lifts the diagonal map $B\to B\times B$ factors as $E\xrightarrow{diag}E\times E\xrightarrow{p\times \id} B\times E$. Recall finally that the cup product in smooth cohomology is defined as the pullback of the exterior product with respect to the diagonal map. We obtain
\begin{equation*}
  \hat{p}_!(p^*x\times y)=(-1)^{kd}x\cup \hat{p}_!(y).
\end{equation*}
   
%
%

\section{Transformations between smooth cohomology}

The construction of smooth cohomology via stratifolds, i.e.~generalized
oriented manifolds, immediately allows to define a lift of the orientation
transformation from a bordism theory which is naturally equipped with an
{$H$}-orientation to the corresponding smooth extensions \textcolor{black}{of
  the present article and of \cite{B-S-S-W}} (provided the characters are chosen appropriately).

As an example, take the canonical orientation from complex bordism to integral homology. As character on complex bordism, use this map composed with the natural map from integral cohomology to cohomology with real coefficients: 
\begin{equation*}
  MU^*(X)\xrightarrow{ori} H^*(X;\bbZ)\xrightarrow{i_*} H^*(X;\bbR).
\end{equation*}

In the stratifold description of integral cohomology and for $X$ an oriented manifold, the transformation sends $[E\to X]\in MU^*(X)$ to $[E\to X]$, where the complex oriented manifold $E$ with proper map to $X$ is interpreted as a stratifold with proper morphism to $X$. A bordism of manifolds over $X$ is also a bordism of stratifolds, so this map is well defined; it is an easy exercise that this indeed describes the natural transformation dual to taking the fundamental class of a stable almost complex manifold. 

We immediately get a smooth lift 
\begin{equation*}
  \hat{MU}^*(X) \to \hat{HS}^*(X)
\end{equation*}
by mapping  the {$ \hat{MU}^*(X)$}-class  $[E,f,\alpha]$ to the  {$ \hat{HS}^*(X)$-class}
$[E,f,\alpha]$.
Obviously this is compatible with the curvature homomorphisms as well as with the passage to the underlying homology theories and the transformation $ori$, and with the action of differential forms. 

\bibliographystyle{plain}




\end{document}